\documentclass{article}

\usepackage{amsfonts}
\usepackage{amsmath}
\usepackage{amssymb}
\usepackage{amsthm}
\usepackage{breakcites}
\usepackage{caption}
\usepackage[utf8]{inputenc}
\usepackage{longtable}
\usepackage{colortbl,xcolor}
\usepackage[a4paper, left=1in, right=1in, bottom=1.2in]{geometry}
\usepackage{graphicx}
\usepackage{hyperref}
\usepackage{natbib}
\usepackage{subfig}
\usepackage{todonotes}

\newcommand{\eps}{\epsilon}
\newcommand{\wnorm}{\tilde{w}}

\DeclareMathOperator{\var}{Var}
\DeclareMathOperator{\cov}{Cov}
\DeclareMathOperator{\cor}{Cor}
\DeclareMathOperator{\expval}{E}
\DeclareMathOperator{\ind}{I}

\hypersetup{
    colorlinks,
    linkcolor={red!70!black},
    citecolor={blue!70!black},
    urlcolor={blue!80!black}
}

\begin{document}

\title{An Analysis of the Johnson-Lindenstrauss Lemma\\ with the Bivariate Gamma Distribution}

\author{Jason\ Bernstein, Alec\ M. Dunton, Benjamin\ W. Priest\\
        Lawrence Livermore National Laboratory\\
        Livermore, CA 94550\\ \\
        \{bernstein8, dunton1, priest2\}@llnl.gov}

\setlength{\parskip}{0.5pc}

\maketitle

\begin{center}
{  \setlength{\baselineskip}{.5\baselineskip}
LLNL-TR-844277
\vskip3.0true in
}
\end{center}

\normalsize

\setlength{\parskip}{0.0pc}

\indent

\vspace{-3in}

\begin{abstract}
Probabilistic proofs of the Johnson-Lindenstrauss lemma imply that random projection can reduce the dimension of a data set and approximately preserve pairwise distances.
If a distance being approximately preserved is called a success, and the complement of this event is called a failure, then such a random projection likely results in no failures.
Assuming a Gaussian random projection, the lemma is proved by showing that the no-failure probability is positive using a combination of Bonferroni's inequality and Markov's inequality.
This paper modifies this proof in two ways to obtain a greater lower bound on the no-failure probability.
First, Bonferroni's inequality is applied to pairs of failures instead of individual failures.
Second, since a pair of projection errors has a bivariate gamma distribution, the probability of a pair of successes is bounded using an inequality from \cite{jensen1969}.
If $n$ is the number of points to be embedded and $\mu$ is the probability of a success, then this leads to an increase in the lower bound on the no-failure probability of $\frac{1}{2}\binom{n}{2}(1-\mu)^2$ if $\binom{n}{2}$ is even and $\frac{1}{2}\left(\binom{n}{2}-1\right)(1-\mu)^2$ if $\binom{n}{2}$ is odd.
For example, if $n=10^5$ points are to be embedded in $k=10^4$ dimensions with a tolerance of $\eps=0.1$, then the improvement in the lower bound is on the order of $10^{-14}$.
We also show that further improvement is possible if the inequality in \cite{jensen1969} extends to three successes, though we do not have a proof of this result.
\end{abstract}

\newpage

{\renewcommand\arraystretch{1.2}
\noindent\begin{longtable*}{@{}l @{\quad:\quad} l@{}}
\multicolumn{2}{@{}l}{\textbf{Notation}} \\
$n\in\mathbb{N}^+$            & Number of points to be embedded \\
$x_j\in\mathbb{R}^d$          & $j$th data point \\
$d\in\mathbb{N}^+$            & Dimension of each $x_j$ \\
$k\in\mathbb{N}^+$            & Dimension of embedding space \\
$R\in\mathbb{R}^{d\times k}$  & Random projection matrix; entries are i.i.d. $N(0,1)$ \\
$w_i=x_j-x_{j'}$              & Difference between $x_j$ and $x_{j'}$ \\
$\wnorm_i=w_i/\|w_i\|_2$      & Unit-vector in direction of $w_i$\\
$\eps\in(0,1)$                & Tolerance on projection error \\
$\Gamma(k/2,2)$               & Gamma distribution with shape parameter $k/2$ and scale parameter 2 \\
$Z_i=R^{\mathsf T}\wnorm_i$   & Random projection, has the $N(0_k,I_k)$ distribution \\
$V_i=\|Z_i\|^2_2$             & Projection error, has the $\Gamma(k/2,2)$ distribution \\
$p(V_i\leq v_i)$              & Cumulative distribution function of a projection error \\
$S_i=\{|V_i-k|\leq k\eps\}$   & The $i$th projection error divided by $k$ is close to one, a `success' \\
$F_i=S_i^c$                   & The complement of $S_i$, a `failure' \\
$\mu=p(S_i)$                  & Success probability, independent of $i$ \\
$p(\cap_i S_i)$               & No-failure probability \\
$p(S_i\cap S_{i'})$           & Joint success probability \\
$h(v_i,v_{i'})$               & Joint probability density function of $V_i$ and $V_{i'}$ \\
$H(v_i,v_{i'})$               & Joint cumulative distribution function of $V_i$ and $V_{i'}$ \\
$\rho_{ii'}^2=(\wnorm_i\cdot\wnorm_{i'})^2$ & Correlation of $V_i$ and $V_{i'}$ \\
$I_k$                         & Identity matrix of dimension $k\times k$ \\
$0_k$                         & Zero vector of length $k$ \\
$\otimes$                     & Kronecker product \\
\multicolumn{2}{@{}l}{\textbf{Subscripts}} \\
$j,j'=1,\ldots,n$             & Indices for individual points \\
$i,i'=1,\ldots,\binom{n}{2}$  & Indices for pairs of points \\
$\ell=1,\ldots,k$             & Index over elements of $Z_i$ \\
\multicolumn{2}{@{}l}{\textbf{Acronyms}} \\
JL                            & Johnson-Lindenstrauss \\
i.i.d.                        & Independent and identically distributed
\end{longtable*}}

\newpage

\section{Introduction}
\label{sec:introduction}

The Johnson-Lindenstrauss (JL) lemma says, informally, that the dimension of a data set can be reduced by projection such that distances between points are almost preserved \citep{jl1984}.
Review papers on the lemma and its extensions include \cite{nelson2020}, \cite{ghojogh2021}, and \cite{freksen2021}, and additional overviews can be found in \cite{vempala2005}, \cite{matouvsek2013}, and \cite{vershynin2018}.
In this paper, the terms success and failure are used to indicate that a distance is preserved, or not, by a projection.
With this terminology, the JL lemma implies that for an arbitrary data set, there is a projection that results in no failures.
Details of the lemma include choosing a sufficiently large embedding dimension, precisely defining what it means for a distance to be preserved, and picking a class of projections.
Due in part to the importance of dimension reduction in data science, research is active on the lemma and related ideas \citep{larsen2017,matouvsek2008}.

A common approach to proving the JL lemma is to show that the no-failure probability is positive for a particular class of random projections \citep{frankl1988,indyk1998,dasgupta2003}.
Note for context that if the projection is random, then the successes and failures are random, and so there is an associated probability of having no failures.
The difficulty with this approach is that the joint distribution of the failures is not available in closed-form, so it is not immediate that the no-failure probability is positive.
However, following \cite{dasgupta2008}, the no-failure probability can be shown to be positive in essentially two steps.
First, Bonferroni's inequality, or union bounding, is used to bound the probability of at least one failure, where the upper bound is the sum of the individual failure probabilities.
Second, an embedding dimension is chosen such that this upper bound is less than one.
This embedding dimension is found using Markov's inequality and the chi-square distribution of the projection errors.
As discussed in detail later, a projection error quantifies how much the random projection changes the distance between two points and determines whether a failure or success occurs.
Applying De Morgan's laws with this choice of embedding dimension then implies that the no-failure probability is positive, which proves the JL lemma.
Overall, this approach bounds the no-failure probability by the sum of the individual failure probabilities and therefore depends on the marginal distribution of the failures.

The contribution of this paper is an improvement to the lower bound on the no-failure probability that appears in the JL lemma.
Our approach is to apply Bonferroni's inequality to pairs of failures, rather than to each failure individually as in \cite{dasgupta2003}.
This leads to computing the joint success probability, or probability that two successes occur, which is a function of the joint distribution of the corresponding projection errors.
The joint distribution is the bivariate gamma distribution studied in \cite{kibble1941}, whose properties have been well-studied and are applicable here.
In particular, an inequality that applies to this distribution from \cite{jensen1969} is used to bound each joint success probability.
Bounding joint success probabilities, rather than probabilities of individual successes, leads to a greater lower bound on the no-failure probability.
The improvement in the lower bound on the no-failure probability is typically small for embedding dimensions and data set sizes that arise in practice.
For example, if $n=10^5$ points are embedded in $k=10^4$ dimensions with tolerance $\eps=0.1$, then the lower bound on the no-failure probability increases by approximately $10^{-14}$.
This improvement was not found to lead to a smaller embedding dimension, which would be desirable from a data compression perspective, but this work may still be a step in that direction.

In the JL lemma, Bonferroni's inequality leads to a conservative bound on the probability of at least one failure, as noted in \cite{li2006a} and \citet[footnote 7]{li2006b}, for example.
The source of this conservativeness is that the inequality is applied to each failure separately, since this ignores information contained in the joint distribution of the failures.
Applying Bonferroni's inequality to individual failures therefore leads to a smaller bound on the no-failure probability than necessary.
Our approach is to apply Bonferroni's inequality to pairs of failures instead, which allows information from the bivariate distributions of the failures to be used.
Combining Bonferroni's inequality in this way with an inequality for bivariate gamma distributions from \cite{jensen1969} leads to a greater lower bound on the no-failure probability.
Sharper Bonferroni inequalities have also been developed that rely on probabilities of pairs of events \citep{hunter1976,worsley1982}, but they are not considered here.
It is also worth noting that much of the statistical work on Bonferroni's inequality, and variations of the inequality, has been related to multiple hypothesis testing.
Somewhat analogously to the situation here, Bonferroni's inequality can lead to conservative hypothesis tests, as noted in \cite{bland1995} and \citet[ch. 10]{wasserman2004}, for example.

As previously noted, the bivariate gamma distribution is the joint distribution of two projection errors resulting from a Gaussian random projection.
For our purposes, the distribution depends on the embedding dimension and a correlation parameter that is the square of the dot product of the points to be embedded.
These dot products can often not be computed exactly due to computational considerations, and so the bound discussed here is independent of the dot products.
The distribution is studied in \cite{wicksell1933}, \cite{kibble1941}, \cite{krishnaiah1961}, \cite{moran1967}, \cite{jensen1969}, and \cite{iliopoulos2005}, and generalized to higher-dimensions in \cite{krishnamoorthy1951}.
\citet[sec. 8.1]{balakrishna2009} and \citet[sec. 48, sec. 2.3]{kotz2004} review the distribution and provide expressions for its density, distribution, and moment generating function.
More generally, \citet[ch. 48]{kotz2004} and \citet[ch. 8]{balakrishna2009} collect these properties for several bivariate distributions that have marginal gamma distributions.
Following \cite{balakrishna2009}, the particular bivariate gamma distribution that occurs here can be called the Kibble-Wicksell bivariate gamma distribution.
We call this distribution the bivariate gamma distribution for simplicity.

The paper is organized as follows.
Section \ref{sec:background} gives background on Gaussian random projections and shows that the bivariate gamma distribution is the joint distribution of two projection errors.
Section \ref{sec:nfp} then uses the inequality from \cite{jensen1969} for bivariate gamma distributions to get an improved lower bound on the no-failure probability.
Section \ref{sec:discussion} provides additional discussion of our analysis and notes limitations of our approach.
Last, Section \ref{sec:conclusion} concludes the paper and considers possible future work.

\section{Random Projections and Projection Errors}
\label{sec:background}

This section reviews linear, Gaussian random projections and shows that pairs of projection errors have the bivariate gamma distribution.
In particular, Section \ref{subsec:randpro} describes linear, Gaussian random projections and shows that the resulting projection errors are marginally gamma distributed.
Section \ref{subsec:bivgamma} describes the bivariate gamma distribution of a pair of projection errors, and Section \ref{subsec:bounding_jsp} reviews an inequality for this distribution that we apply to the JL lemma.
While the material may not be new to some readers, reviewing this section may still be useful since it introduces notation and terminology used in later sections.

\subsection{Linear, Gaussian Random Projection}
\label{subsec:randpro}

This paper considers random projections that can be represented as a matrix of independent random variables, each having the $N(0,1)$ distribution.
If the matrix is $R\in\mathbb{R}^{d\times k}$ and $w\in\mathbb{R}^d$, then $R^{\mathsf T}w\in\mathbb{R}^k$ is the random projection of $w$ and is a Gaussian random variable by linearity.
The point $R^{\mathsf T}w$ is referred to as the embedding of $w$ and $k$ is referred to as the embedding dimension.
Random projections have several favorable properties that make them useful for dimension reduction, in which case $k<<d$.
For example, they do not depend on the data directly and they are not likely to significantly change distances between points if the embedding dimension is sufficiently large.
More detailed overviews of random projections can be found in \cite{nelson2020} and \cite{ghojogh2021}, for example, and additional references are given in Section \ref{sec:introduction}.

For notation, let $x_j$ and $x_{j'}$ be two of the $n$ points to be embedded, and let $i=i(j,j')$ be an index over pairs of points.
The difference between the $i$th pair of points is denoted $w_i=x_j-x_{j'}$, so that $\wnorm_i=w_i/\|w_i\|_2$ is the difference between the pair of points scaled to have length one.
The $i$th random projection is defined as
\begin{equation}
Z_i=R^{\mathsf T}\wnorm_i,
\end{equation}
and is Gaussian distributed by linearity,
\begin{equation}
Z_i\sim N(0_k,I_k).
\end{equation}
Since the elements of the random projection vector are independent, $N(0,1)$ random variables, it follows that the $i$th projection error, defined as
\begin{equation}
\label{Vi}
V_i=\|Z_i\|^2_2,
\end{equation}
has the chi-square distribution with $k$ degrees of freedom.
Equivalently, $V_i$ has the gamma distribution with shape parameter $k/2$ and scale parameter 2, denoted $V_i\sim\Gamma(k/2,2)$.
The random variable $V_i$ is called a projection error because it quantifies how much the random projection changes the distance between two points.
Note that the mean of a projection error is $\expval[V_i]=k$ and the variance of a projection error is $\var[V_i]=2k$.
It is common to see the projection matrix defined instead as $R/\sqrt{k}$, which leads to the projection errors having mean one and variance $2/k$.
The different scalings lead to different interpretations of the projection errors but do not change our results.
The scaling employed here is used since it simplifies the modeling of the projection errors with the bivariate gamma distribution.

The important point is that if $V_i/k$ is close to one, then the random projection, scaled by $1/\sqrt{k}$, approximately preserved the distance between the points $x_j$ and $x_{j'}$.
Hence, if we define a success as the event that a projection error is close to its mean, then ideally the success probability is high.
The notion of a projection error being close to its mean is formalized by defining the $i$th success as
\begin{equation}
S_i=\{V_i\in[k(1-\eps),k(1+\eps)]\},
\end{equation}
where $\eps\in(0,1)$ is called the tolerance since it determines the limits on the allowed projection error.
The probability of a success, $p(S_i)$, depends on the embedding dimension and the tolerance, but not the points to be embedded or index $i$, and so is a constant,
\begin{equation}
\label{mu}
\mu=p(V_i\leq k(1+\eps)) - p(V_i\leq k(1-\eps)).
\end{equation}
By the weak law of large numbers, the success probability converges to one as $k\rightarrow\infty$.
A failure is defined as the complement of a success.
The $i$th failure is denoted $F_i=S_i^c$ and the probability of this outcome is $p(F_i)=1-\mu$ by \eqref{mu}.

\begin{figure}[ht]
\includegraphics[scale=0.5]{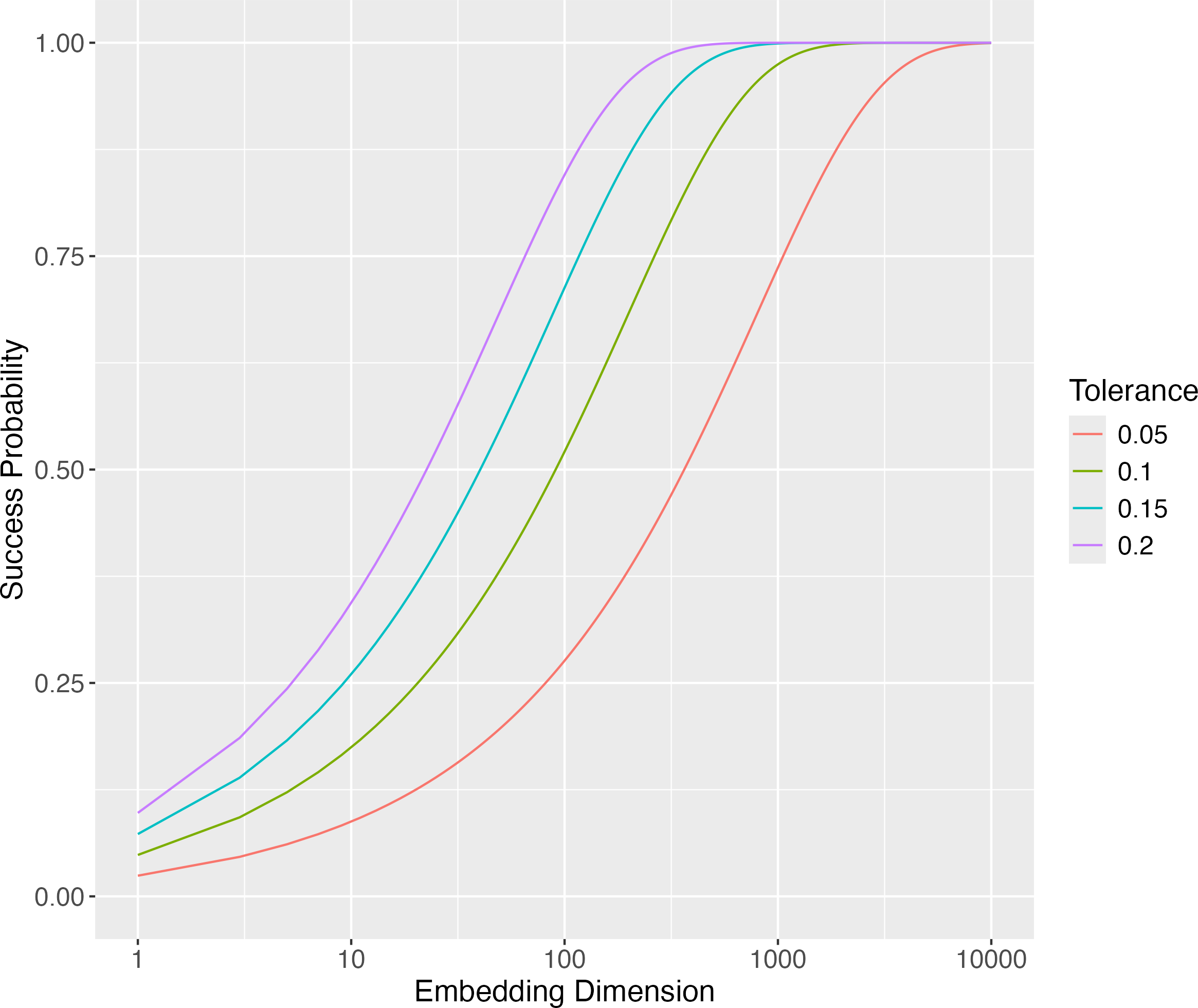}
\centering
\caption{The success probability for different embedding dimensions and tolerances.}
\label{fig:success_probability}
\end{figure}

Figure~\ref{fig:success_probability} shows how the success probability changes with the embedding dimension and tolerance.
We see that the success probability converges to one as $k\rightarrow\infty$ for fixed $\eps$ and is positive for all possible values of the embedding dimension.
The probabilities are computed using the $\Gamma(k/2,2)$ distribution function.

\subsection{Bivariate Gamma Distribution}
\label{subsec:bivgamma}

The approach taken here depends on bounding probabilities of the form $p(S_i\cap S_{i'})$, which is a function of the joint distribution of $V_i$ and $V_{i'}$.
We refer to these probabilities as joint success probabilities, in contrast to the marginal success probability $\mu=p(S_i)$ defined in Section \ref{subsec:randpro}.
A brief review of the main points of this section is now provided to help guide the discussion:
\begin{enumerate}
\item The joint distribution of two projection errors, $V_i$ and $V_{i'}$, is the bivariate gamma distribution described in \cite{kibble1941}.
This distribution depends on $k$ and the correlation of $V_i$ and $V_{i'}$, which is the squared-dot product of $\wnorm_i$ and $\wnorm_{i'}$, or $(\wnorm_i\cdot \wnorm_{i'})^2$;
\item Joint success probabilities can thus be computed if $w_i$ and $w_{i'}$ are known, but this may not be possible if the data has not been collected yet or $n$ is large, for example;
\item Instead, the joint success probabilities are bounded below using an inequality from \cite{jensen1969}.
The inequality for our purposes and in our notation is $p(S_i\cap S_{i'})\geq\mu^2$.
\end{enumerate}
The remainder of this section provides additional details on these points.

To identify the joint distribution of the projection errors $V_i$ and $V_{i'}$, it is helpful to start with the joint distribution of the random projections $Z_i$ and $Z_{i'}$.
By linearity and independence of the elements of $R$, the two random projections are jointly Gaussian distributed,
\begin{equation}
\label{joint_normal}
\begin{pmatrix}
Z_i \\ Z_{i'}
\end{pmatrix}
\sim N\left(
\begin{pmatrix}
0_k \\ 0_k
\end{pmatrix},
\begin{pmatrix}
I_k & \rho_{ii'} I_k \\
\rho_{ii'} I_k & I_k
\end{pmatrix}
\right),
\end{equation}
where $I_k$ is the identity matrix of dimension $k\times k$ and the correlation parameter in the covariance matrix is
\begin{equation}
\label{rho}
\rho_{ii'}=\wnorm_i\cdot\wnorm_{i'}.
\end{equation}
The correlation between $Z_i$ and $Z_{i'}$ is due to each being a function of the same random matrix $R$.
The covariance matrix in \eqref{joint_normal} can also be written as
\begin{equation}
\Sigma_{ii'}\otimes I_k,
\end{equation}
where
\begin{equation}
\label{Sigma_ii}
\Sigma_{ii'}=
\begin{pmatrix}
1 & \rho_{ii'} \\
\rho_{ii'} & 1
\end{pmatrix}
\end{equation}
is a correlation matrix and $\otimes$ is the Kronecker product.
Hence, the random projections $Z_i$ and $Z_{i'}$ have the same distribution as $k$ independent samples from a bivariate normal distribution with covariance matrix $\Sigma_{ii'}$.
The full joint distribution of all $\binom{n}{2}$ random projections is given in Appendix~\ref{subsec:joint_dist} for completeness.
Joint distributions of random projections are also used for estimation in \cite{li2006a}, \cite{li2006b}, and \cite{kang2021}, for example.

The Gaussian distribution of the random projections in \eqref{joint_normal} determines the joint distribution of the projection errors $V_i$ and $V_{i'}$.
Specifically, the joint distribution of $V_i$ and $V_{i'}$ is the bivariate gamma distribution discussed in \cite{kibble1941} and reviewed in \citet[sec. 8.2]{balakrishna2009}.
Furthermore, since $V_i$ and $V_{i'}$ are each the sum of $k$ independent, squared, standard normal random variables, they are both marginally $\Gamma(k/2,2)$ distributed.
The projection errors are independent if and only if $\rho_{ii'}=0$ since Gaussian random variables are independent if and only if they are uncorrelated.
This correlated gamma distribution characterization of $V_i$ and $V_{i'}$ leads to their joint distribution being called the bivariate gamma distribution.

The correlation of $V_i$ and $V_{i'}$ is computed from the moment generating function of the bivariate gamma distribution with parameters $k$ and $\rho_{ii'}$.
A derivation of the moment generating function is given in \cite{kibble1941} and \cite{krishnamoorthy1951}, and Appendix \ref{subsec:charfun} for completeness.
To derive the correlation of $V_i$ and $V_{i'}$, note first from \eqref{joint_normal} that
\begin{equation}
\label{gaussian_V}
\begin{pmatrix}
V_i \\ V_{i'}
\end{pmatrix}
\stackrel{d}{=}\sum_{\ell=1}^k
\begin{pmatrix}
Z_{i\ell}^2 \\ Z_{i'\ell}^2
\end{pmatrix},
\end{equation}
where $Z_{i\ell}$ and $Z_{i'\ell}$ are the $\ell$th elements of $Z_i$ and $Z_{i'}$, respectively, and $\cov[Z_{i\ell},Z_{i'\ell}]=\rho_{ii'}$.
It is shown in Appendix \ref{subsec:charfun} that $\cov[Z_{i\ell}^2,Z_{i'\ell}^2]=2\rho_{ii'}^2$, which with \eqref{gaussian_V} implies that the covariance of $V_i$ and $V_{i'}$ is
\begin{equation}
\label{cov_proj_errors}
\cov[V_i,V_{i'}]=2k\rho_{ii'}^2.
\end{equation}
The correlation of $V_i$ and $V_{i'}$ is thus
\begin{align}
\cor[V_{i},V_{i'}]
&=\frac{\cov[V_i,V_{i'}]}{(\var[V_i]\var[V_{i'}])^{1/2}}\\
&=\rho_{ii'}^2
\end{align}
since $\var[V_i]=\var[V_{i'}]=2k$.
It follows that $V_i$ and $V_{i'}$ are positively correlated if $\rho_{ii'}\neq0$, or equivalently if $w_i$ and $w_{i'}$ are not orthogonal.

The bivariate gamma density function of the projection errors is derived in \cite{kibble1941}.
Specifically, \citet[eqn. 12]{kibble1941} gives the density of $V_i/2$ and $V_{i'}/2$, so a scaling by two is applied to get the joint density of $V_i$ and $V_{i'}$ as
\begin{equation}
\label{pdf_kibble}
h(v_i,v_{i'})=\frac{1}{4\Gamma(k/2)(1-\rho_{ii'}^2)}\left(\frac{v_iv_{i'}}{4\rho_{ii'}^2}\right)^{(k/2-1)/2}\exp\left(-\frac{v_i+v_{i'}}{2(1-\rho_{ii'}^2)}\right)I_{k/2-1}\left(\frac{\sqrt{v_iv_{i'}\rho_{ii'}^2}}{1-\rho_{ii'}^2}\right),
\end{equation}
where $I_{k/2-1}(\cdot)$ is the modified Bessel function of the first kind and order $k/2-1$;
see also the last equation in \citet[sec. 5.4]{kibble1941} for a similar expression.
The scaling by one-half appears to be a difference between the bivariate gamma distribution and random projection literature, which we follow the latter on and do not scale.
If $\rho_{ii'}=0$, then $V_i$ and $V_{i'}$ are independent and the joint density function is the product of the marginal $\Gamma(k/2,2)$ density functions of $V_i$ and $V_{i'}$.
This can be obtained from the joint density \eqref{pdf_kibble} using the identity
\begin{equation}
I_{k/2-1}\left(\frac{\sqrt{v_iv_{i'}\rho_{ii'}^2}}{1-\rho_{ii'}^2}\right)\sim \frac{1}{\Gamma(k/2)}\left(\frac{1}{2}\frac{\sqrt{v_iv_{i'}\rho_{ii'}^2}}{1-\rho^2_{ii'}}\right)^{(k/2-1)}
\end{equation}
as $\rho_{ii'}\rightarrow0$ \citep[eqn. 9.6.7]{abramowitz1964}.
\cite{iliopoulos2005} gives expressions for the conditional means and variances of the distribution.

\begin{figure}[ht]
\includegraphics[scale=0.5]{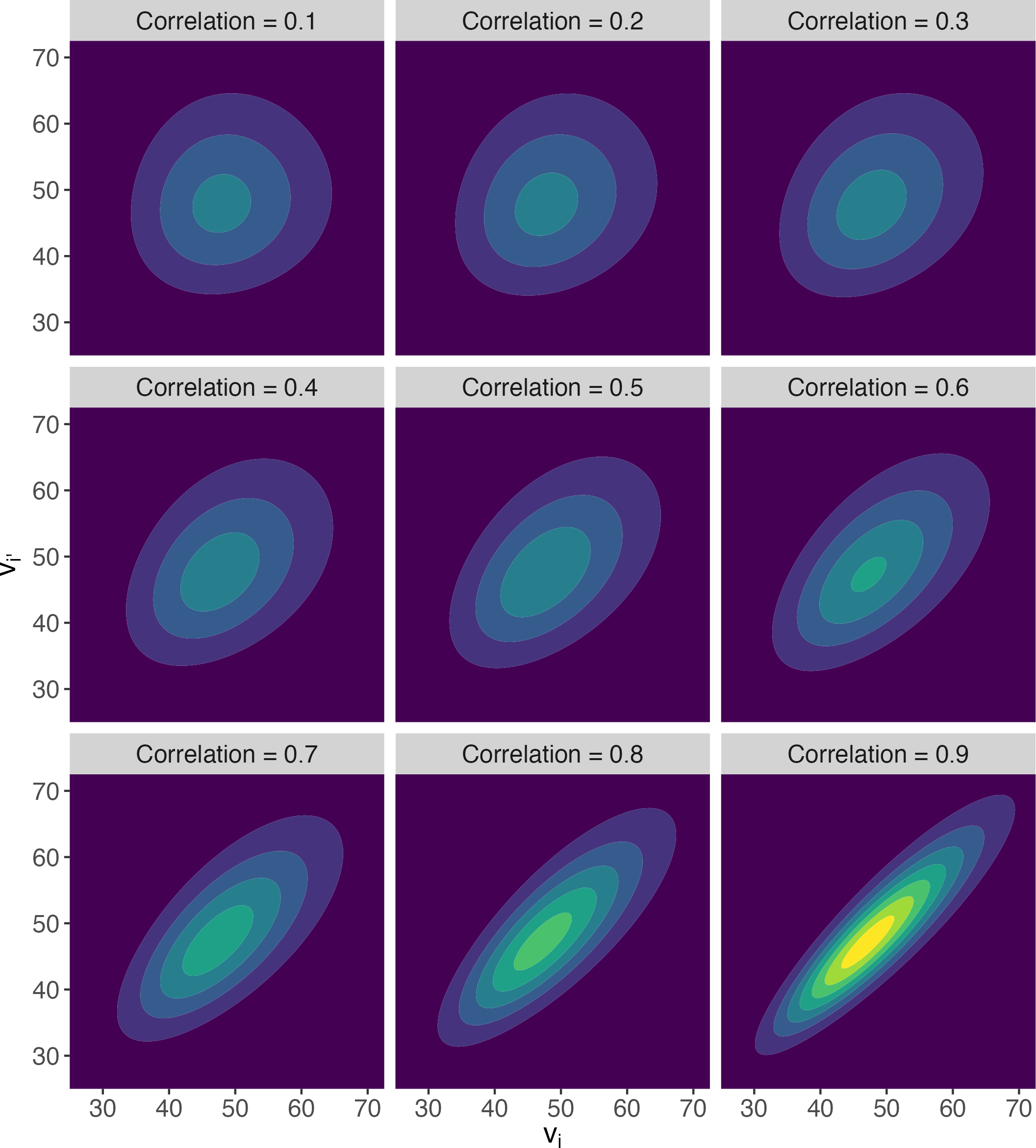}
\centering
\caption{The bivariate gamma density function for different values of the correlation parameter, $\rho_{ii'}^2$, and an embedding dimension of $k=50$.  The marginal densities are $\Gamma(25,2)$ for all $\rho_{ii'}^2$.}
\label{fig:bivgamma}
\end{figure}

Figure \ref{fig:bivgamma} shows the bivariate gamma density function as the correlation, $\rho_{ii'}^2$, is varied from 0.1 to 0.9 in increments of 0.1.
The density is symmetric in $v_i$ and $v_{i'}$ and becomes increasingly concentrated along the $v_i=v_{i'}$ axis as $\rho_{ii'}^2$ increases.
Recall as well that the marginal expectations and variances are equal to $k$ and $2k$ since the projection errors have the $\Gamma(k/2,2)$ distribution.

\begin{figure}[ht]
\includegraphics[scale=0.45]{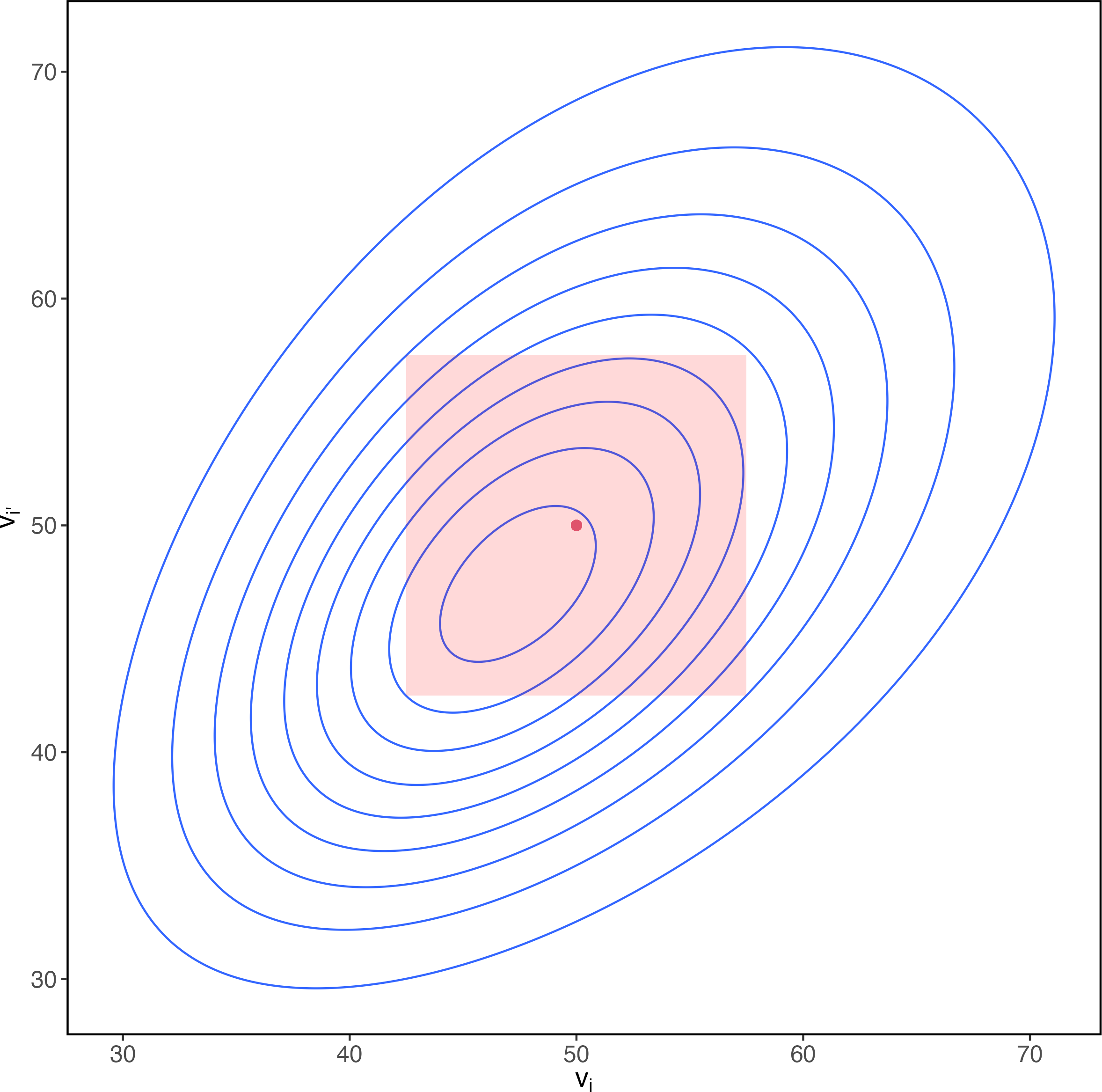}
\centering
\caption{A bivariate gamma distribution with $k=50$ and $\rho_{ii'}^2=0.5$.  The red shaded area indicates the integration domain for computing the joint success probability.  In each coordinate, the integration limits are $(k(1-\eps),k(1+\eps))$, where $\eps=0.15$.  The red dot is at $(k,k)$ and indicates the marginal means.}
\label{fig:integration_domain}
\end{figure}

Figure \ref{fig:integration_domain} shows the integration domain for computing the joint success probability.
The domain is square since $S_i\cap S_{i'}$ is equivalent to $V_i\in[k(1-\eps),k(1+\eps)]$ and $V_{i'}\in[k(1-\eps),k(1+\eps)]$, and so the joint success probability is the integral
\begin{equation}
\label{joint_succ_prob}
p(S_i\cap S_{i'})=\int_{k(1-\eps)}^{k(1+\eps)}\int_{k(1-\eps)}^{k(1+\eps)}h(v_i,v_{i'})\,dv_i\,dv_{i'}.
\end{equation}
Alternatively, if the joint cumulative distribution function of $V_i$ and $V_{i'}$ is denoted
\begin{equation}
H(v_i,v_{i'})=p(V_i\leq v_i, V_{i'}\leq v_{i'}),
\end{equation}
then the joint success probability can be written as
\begin{equation}
\label{joint_cdf}
p(S_i\cap S_{i'})=H(k(1+\eps),k(1+\eps))-2H(k(1+\eps),k(1-\eps))+H(k(1-\eps),k(1-\eps)),
\end{equation}
where the middle term follows from the symmetry of the distribution function, $H(v_i,v_{i'})=H(v_{i'},v_i)$.
Formulas for $H(v_i,v_{i'})$ are given in \citet[sec. 8.2.2]{balakrishna2009} and the references therein.
The joint success probability can thus, in principle, be computed by numerically integrating the density in \eqref{pdf_kibble} or evaluating \eqref{joint_cdf}.
However, numerically integrating the density function \eqref{pdf_kibble} may be difficult when $k$ is large, and so an approximation of the distribution function is needed in this case.
To this end, the central limit theorem implies that the distribution of $V_i$ and $V_{i'}$, suitably scaled, converges to a bivariate Gaussian distribution as $k\rightarrow\infty$.
The bivariate gamma distribution $H(v_i,v_{i'})$ can therefore be quickly approximated using a bivariate normal approximation for sufficiently large $k$.
Details on this approach are provided in Appendix \ref{subsec:clt}.

\subsection{Bounding the Joint Success Probability}
\label{subsec:bounding_jsp}

An inequality for bivariate gamma distributions from \cite{jensen1969} provides a lower bound on the joint success probabilities.
The inequality applies here since, by \eqref{joint_normal}, $Z_i$ and $Z_{i'}$ are jointly Gaussian distributed and each has expected value $0_k$ and a diagonal covariance matrix.
Under these conditions, \cite{jensen1969} showed that if $a\geq b\geq0$, then
\begin{equation}
p(V_i\in[a,b], V_{i'}\in[a,b])\geq p(V_i\in[a,b]) p(V_{i'}\in[a,b]),
\end{equation}
where $V_i=\|Z_i\|^2_2$ and $V_{i'}=\|Z_{i'}\|^2_2$.
Taking $a=k(1-\eps)$ and $b=k(1+\eps)$ in this inequality, it follows that the joint success probabilities are bounded below by $\mu^2$, or
\begin{equation}
\label{jensen}
p(S_i\cap S_{i'})\geq\mu^2,
\end{equation}
with equality if and only if $\rho_{ii'}=0$.
The equality case is verified by noting that $V_i$ and $V_{i'}$ are independent if and only if $\rho_{ii'}=0$, in which case $p(S_i\cap S_{i'})=p(S_i)p(S_{i'})=\mu^2$.
Furthermore, the joint success probability is always positive since $\mu>0$ for all tolerances and embedding dimensions, which implies that two successes are never mutually exclusive.
It should also be noted that the inequality in \cite{jensen1969} applies to a broader range of bivariate gamma distributions than the particular distribution considered here.
Specifically, the inequality applies even if the covariance matrix of $Z_i$ and $Z_{i'}$ does not have the form $\rho_{ii'}I_k$, though this full generality is not needed here.

The inequality on joint success probabilities \eqref{jensen} is equivalent to
\begin{equation}
\label{condl_succ_prob}
p(S_{i'}|S_i)\geq\mu
\end{equation}
since $p(S_i\cap S_{i'})=p(S_{i'}|S_i)\mu$ by Bayes' rule, meaning that the conditional success probability $p(S_{i'}|S_i)$ can be greater than the success probability.
In words, this inequality says that knowledge of one success may increase the probability of a second success.
This suggests another interpretation of \eqref{jensen}, which is that the success-indicator random variables, $\ind[S_i]$ and $\ind[S_{i'}]$, are non-negatively correlated.

An upper bound on the joint success probability is also available.
Since $\{S_i\cap S_{i'}\}\subseteq S_i$ and $p(S_i)=\mu$, it follows that $\mu$ is an upper bound on $p(S_i\cap S_{i'})$, or
\begin{equation}
\label{mu_inequality}
p(S_i\cap S_{i'})\leq\mu.
\end{equation}
Note that if $\rho_{ii'}^2=1$, then $\wnorm_i=\wnorm_{i'}$ or $\wnorm_i=-\wnorm_{i'}$, implying $V_i=V_{i'}$ and thus $p(S_i\cap S_{i'})=p(S_i)=\mu$.
Hence, equality occurs in \eqref{mu_inequality} if $\rho_{ii'}^2=1$, or $w_i$ and $w_{i'}$ point in the same or opposite directions.
In summary, $\mu^2\leq p(S_i\cap S_{i'})\leq\mu$, with the lower and upper bounds obtained when $\rho_{ii'}=0$ and $\rho_{ii'}^2=1$, respectively.

\begin{figure}[ht]
\includegraphics[scale=0.5]{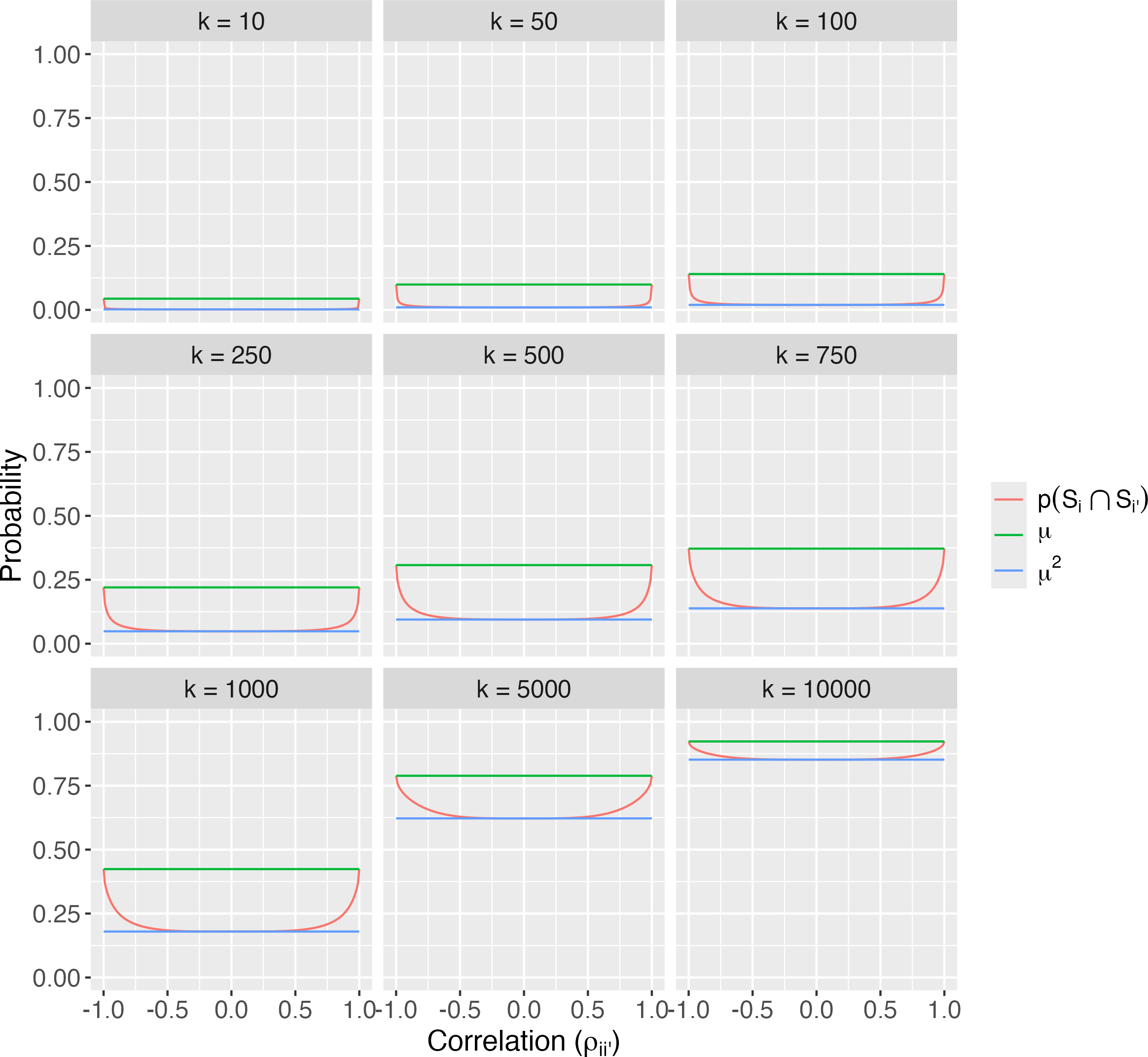}
\centering
\caption{The joint success probability as a function of the embedding dimension and the correlation $\rho_{ii'}$.  The probability equals $\mu$ when $\rho_{ii'}^2=1$ and equals $\mu^2$ when $\rho_{ii'}=0$.  The tolerance here is $\eps=0.025$.}
\label{fig:jsp}
\end{figure}

Figure \ref{fig:jsp} shows how the joint success probability changes as a function of the correlation, $\rho_{ii'}$, and embedding dimension, $k$.
For all $k$, the joint success probability equals $\mu^2$ and $\mu$ when $\rho_{ii'}=0$ and $\rho_{ii'}^2=1$, respectively, and these are the lower and upper bounds on the probability.
The bounds vary with the embedding dimension, and for example are $(0.18,0.42)$ when $k=10^3$ and $(0.85,0.92)$ when $k=10^4$.
To create this plot, a Gaussian approximation of the joint success probability was made when $k>300$ to avoid issues with numerically integrating the density in \eqref{pdf_kibble}.
The appendix provides details on this approximation, which is justified by the central limit theorem as noted at the end of Section \ref{subsec:bivgamma}.

We continue with reviewing the JL lemma, which shows that the no-failure probability is positive for a particular embedding dimension.
The inequality in \eqref{jensen} is then used to improve the lower bound on the no-failure probability.

\section{Bounding the No-Failure Probability}
\label{sec:nfp}

The JL lemma gives a lower bound on the embedding dimension, as a function of $n$ and $\eps$, that ensures the no-failure probability is positive,
\begin{equation}
p\left(\cap_{i=1}^{\binom{n}{2}}S_i\right)>0.
\end{equation}
The no-failure probability can quickly be shown to be positive if the successes are mutually independent or if all $\rho_{ii'}^2=1$.
That is, if the successes are mutually independent, or all $\rho_{ii'}=0$, then the no-failure probability is $\mu^{\binom{n}{2}}$ and positive since $\mu>0$ for all $k$ and $\eps$.
Alternatively, if all $\rho_{ii'}^2=1$, then the no-failure probability equals $\mu$ and is again positive for all values of $k$ and $\eps$.
However, if some $\rho_{ii'}^2\notin\{0,1\}$, then the distribution of the number of successes is not available in closed-form and the no-failure probability is not clearly positive for a given $k$ and $\eps$;
recall from Section \ref{subsec:bivgamma} that the random projections are correlated due to the random matrix $R$, which implies that the successes are dependent as well.
An approach to showing that the no-failure probability is positive is therefore needed that allows for dependent successes.

In this section, we first review the approach in \cite{dasgupta2003} that uses Bonferroni's inequality to get a lower bound on the no-failure probability.
The approach is then adapted to find a greater lower bound on the no-failure probability using the lower bound on the joint success probability given in \eqref{jensen}.
The new lower bound is greater than the standard lower bound for a fixed embedding dimension, though the difference depends on the embedding dimension as discussed later in this section.

We begin by reviewing the proof of the JL lemma in \cite{dasgupta2003} to show that the no-failure probability is positive if the embedding dimension is sufficiently large.
The first step is to bound the probability of one or more failures using Bonferroni's inequality,
\begin{align}
\label{bonf}
p\left(\cup_{i=1}^{\binom{n}{2}}F_i\right)
&\leq \sum_{i=1}^{\binom{n}{2}}p(F_i)\\
&=\binom{n}{2}(1-\mu),
\end{align}
which shows that the probability of at least one failure is bounded above by the expected number of failures.
Taking complements and using De Morgan's laws, it follows that the no-failure probability satisfies
\begin{equation}
\label{mlb}
p\left(\cap_{i=1}^{\binom{n}{2}}S_i\right)\geq1-\binom{n}{2}(1-\mu).
\end{equation}
The inequality in \eqref{mlb} can be written in an equivalent form that is commonly-used,
\begin{equation}
p\left(\cap_{i=1}^{\binom{n}{2}}S_i\right)\geq\sum_{i=1}^{\binom{n}{2}}p(S_i)-\left(\binom{n}{2}-1\right),
\end{equation}
where the substitution $p(S_i)=\mu$ is made to complete the equivalence \citep[eqn. 1.2.10]{casella2021}.
We refer to the lower bound in \eqref{mlb} as the marginal lower bound since it is derived with the marginal distribution of the projection errors.
Furthermore, since $\mu\rightarrow1$ as $k\rightarrow\infty$ for fixed $\eps$, an embedding dimension exists such that the marginal lower bound is positive;
see \citet[eqn. 2.1]{dasgupta2003} for a lower bound on the embedding dimension that implies the no-failure probability is positive.
Note as well that since the marginal lower bound is obtained by applying Bonferroni's inequality to each failure separately, any possible dependency between failures is ignored.

A greater lower bound on the no-failure probability is found by repeating the steps above, but applying Bonferroni's inequality to pairs of failures instead of each failure individually.
The inequality \eqref{jensen} from \cite{jensen1969} is then applied to each of the joint success probabilities appearing in the resulting inequality.
Bounding the probability of pairs of failures instead of single failures allows the joint distribution of the projection errors to be used in bounding the no-failure probability.
The number of failures to consider, $\binom{n}{2}$, is assumed to be even initially to simplify the notation and derivation of the lower bound.
Proceeding as for the marginal lower bound, the first step in deriving the improved lower bound is to get an upper bound on the probability of at least one failure,
\begin{align}
p(\cup_{i=1}^{\binom{n}{2}}F_i)
&=p(\cup_{i=1}^{\frac{1}{2}\binom{n}{2}}\{F_{2i-1}\cup F_{2i}\})\\
\label{bonf2}
&\leq\sum_{i=1}^{\frac{1}{2}\binom{n}{2}}p(F_{2i-1}\cup F_{2i})\\
\label{DeMorgan2}
&=\sum_{i=1}^{\frac{1}{2}\binom{n}{2}}\left(1-p(S_{2i-1}\cap S_{2i})\right)\\
\label{jensen2}
&\leq\sum_{i=1}^{\frac{1}{2}\binom{n}{2}}(1-\mu^2)\\
\label{cnfp}
&=\frac{1}{2}\binom{n}{2}(1-\mu^2),
\end{align}
where \eqref{bonf2} follows from Bonferroni's inequality, \eqref{DeMorgan2} follows from De Morgan's laws, and \eqref{jensen2} follows from \eqref{jensen}.
Taking the complement of the first expression and applying De Morgan's laws again gives
\begin{equation}
\label{blb}
p\left(\cap_{i=1}^{\binom{n}{2}}S_i\right)\geq1-\frac{1}{2}\binom{n}{2}(1-\mu^2).
\end{equation}
The lower bound in \eqref{blb} is referred to as the bivariate lower bound since it is derived with the bivariate gamma distribution of the projection errors.
As with the marginal lower bound, an embedding dimension that makes the left side of \eqref{blb} positive can be found since $\mu\rightarrow1$ as $k\rightarrow\infty$.
Note as well that the bivariate lower bound is a worst-case bound on the no-failure probability in that it corresponds to the case where all $w_i\cdot w_{i'}=0$.
If instead at least one $w_i\cdot w_{i'}\neq0$, then $p(S_i\cap S_{i'})>\mu^2$ since $\rho_{ii'}^2>0$ and the inequality in \eqref{blb} is strict.

The bivariate lower bound is strictly greater than the marginal lower bound for all embedding dimensions and tolerances.
This is verified by subtracting the marginal lower bound in \eqref{mlb} from the bivariate lower bound in \eqref{blb} and noting that the difference,
\begin{equation}
\label{diff_bounds}
\Delta(n,\mu)=\frac{1}{2}\binom{n}{2}(1-\mu)^2,
\end{equation}
is positive for all $\mu$, though it does converge to zero as the embedding dimension increases since $\mu\rightarrow1$ as $k\rightarrow\infty$.
In general, $k$ will be taken to be large enough that the marginal lower bound is positive, ensuring that the no-failure probability is positive.
This makes the success probability $\mu$ close to one, which leads to a small difference between the bivariate and marginal lower bounds.
For example, the difference is approximately $\Delta(n,\mu)\approx10^{-14}$ when $n=10^5$, $k=10^4$, and $\eps=0.1$, as noted in Section \ref{sec:introduction}.

If the number of failures is odd, then one failure cannot be paired off with another failure and the upper bound \eqref{cnfp} and difference \eqref{diff_bounds} are modified accordingly.
Specifically, if $\binom{n}{2}$ is odd, then the upper bound on the no-failure probability given in \eqref{cnfp} is $\frac{1}{2}\left(\binom{n}{2}-1\right)(1-\mu^2)+(1-\mu)$.
In this expression, the first summand corresponds to $\binom{n}{2}-1$ pairs of failures and the second summand corresponds to the failure that cannot be paired off with another failure.
The difference between the bivariate lower bound and the marginal lower bound is given by \eqref{diff_bounds} with the $\binom{n}{2}$ replaced by $\binom{n}{2}-1$, which is positive.
Going forward, the number of failures to consider is assumed to be even for simplicity.

A greater lower bound on the no-failure probability than the bivariate lower bound is available if the points to be embedded have a certain geometry.
In particular, since each $p(S_{2i-1}\cap S_{2i})\leq\mu$ by \eqref{mu_inequality}, it follows by summing over joint success probabilities that
\begin{equation}
\label{upbnd}
1 - \frac{1}{2}\binom{n}{2}(1-\mu)\geq 1 - \sum_{i=1}^{\frac{1}{2}\binom{n}{2}}\left(1-p(S_{2i-1}\cap S_{2i})\right),
\end{equation}
where equality holds if and only if $\rho_{2i-1,2i}^2=1$ for all $i$ and the lower bound is from \eqref{DeMorgan2}.
Hence, if the points can be indexed such that $\rho_{{2i-1},{2i}}^2=1$ for all $i$, then the no-failure probability is greater than or equal to the left side of \eqref{upbnd}.
However, there is no guarantee in general that the left side of \eqref{upbnd} is a lower bound on the no-failure probability, whereas the bivariate lower bound always holds.

To summarize this section, the no-failure probability can be shown to be positive using the marginal or bivariate distributions of the projection errors.
The two lower bounds on the no-failure probability are derived similarly using Bonferroni's inequality and De Morgan's laws.
However, the bivariate lower bound is greater than the marginal lower bound by an amount that depends on the number of points to be embedded and the success probability.
The difference between the bounds is due to a lower bound on the joint success probabilities that is obtained only when the projection errors are independent.
If the projection errors are correlated, then the bivariate lower bound is conservative by an amount that increases with the magnitudes of the correlations.
Hence, accounting for possible dependency or correlation between successes leads to an improved lower bound on the no-failure probability.

\section{Discussion}
\label{sec:discussion}

This section discusses several additional points related to our approach and main result.
First, Section \ref{subsec:generalizing_the_derivation} considers whether the approach can be extended by applying Bonferroni's inequality to more than two successes at a time.
Section \ref{subsec:selection_of_embedding_dimension} then discusses whether the improved lower bound on the no-failure probability affects selection of the embedding dimension.
Section \ref{subsec:data_bound} gives a data-dependent lower bound on the no-failure probability that improves on the bivariate lower bound.
Unfortunately the data-dependent lower bound cannot be computed when $n$ is large since it is a summation over $\binom{n}{2}$ terms, though it could potentially be estimated.
Last, Section \ref{subsec:limitations} notes limitations of our approach.

\subsection{Generalizing the Derivation}
\label{subsec:generalizing_the_derivation}

A natural follow-up question to this work is whether the derivation can be generalized by applying Bonferroni's inequality to three or more successes at a time.
The answer depends on whether the inequality in \cite{jensen1969} extends to higher-dimensional gamma distributions, which we do not know.
However, the inequality for three successes,
\begin{equation}
\label{jensen3}
p(S_{3i-2}\cap S_{3i-1}\cap S_{3i})\geq \mu^3.
\end{equation}
can be shown to hold in a few simple cases, where we take $i=1$ for simplicity:
\begin{enumerate}
\item If $S_1,S_2,S_3$ are independent, then $\rho_{12}=\rho_{23}=\rho_{31}=0$ and $p(S_1\cap S_2\cap S_3)=\mu^3$;
\item If $\rho_{12}^2=\rho_{23}^2=\rho_{31}^2=1$, then $p(S_1\cap S_2\cap S_3)=\mu$ and \eqref{jensen3} holds since $\mu>\mu^3$;
\item If $S_1$ and $S_2$ are conditionally independent given $S_3$, then $p(S_1\cap S_2|S_3)=p(S_1|S_3)p(S_2|S_3)$ and
\begin{align}
p(S_1\cap S_2\cap S_3)
&=p(S_1\cap S_2|S_3)p(S_3)\\
&=p(S_1|S_3)p(S_2|S_3)p(S_3)\\
&\geq \mu^3,
\end{align}
where the second equality follows from conditional independence and the inequality follows from \eqref{condl_succ_prob} and $p(S_3)=\mu$.
The condition for $S_1$ and $S_2$ to be conditionally independent is derived in Section \ref{subsec:joint_dist} as
\begin{equation}
\rho_{12\cdot3}=0,
\end{equation}
where $\rho_{12\cdot3}$ is the partial correlation coefficient of $Z_1$ and $Z_2$ given $Z_3$ and defined in \eqref{partial_corr_coef}.
\end{enumerate}
It seems plausible that \eqref{jensen3} holds in general since, by \eqref{jensen}, it is equivalent to
\begin{equation}
p(S_3|S_1\cap S_2)\geq\mu,
\end{equation}
which says that the probability of a success conditional on two successes is at least as great as the unconditional probability of a success.
However, we do not have a proof of \eqref{jensen3} for arbitrary correlations.

We can quantify how the lower bound on the no-failure probability changes if the inequality in \cite{jensen1969} extends to any three successes, or if \eqref{jensen3} always holds.
Redoing the derivation of the bivariate lower bound with \eqref{jensen3} yields the following lower bound on the no-failure probability,
\begin{equation}
p(\cap_{i=1}^{\binom{n}{2}} S_i)\geq 1 - \frac{1}{3}\binom{n}{2}(1-\mu^3).
\end{equation}
The difference between this lower bound and the bivariate lower bound is
\begin{equation}
\binom{n}{2}\left(\frac{1}{3}\mu^3-\frac{1}{2}\mu^2+\frac{1}{6}\right),
\end{equation}
which is positive, decreasing over $\mu\in[0,1]$, and converges to zero as $\mu\rightarrow1$.
Hence, if the inequality in \cite{jensen1969} generalized to \eqref{jensen3}, then the lower bound on the no-failure probability would be improved relative to the bivariate lower bound.
We leave the tasks of proving or disproving \eqref{jensen3} and determining if there is further value in extending this argument to more than three successes as future work.

\subsection{Selection of Embedding Dimension}
\label{subsec:selection_of_embedding_dimension}

\begin{figure}[ht]
\includegraphics[scale=0.5]{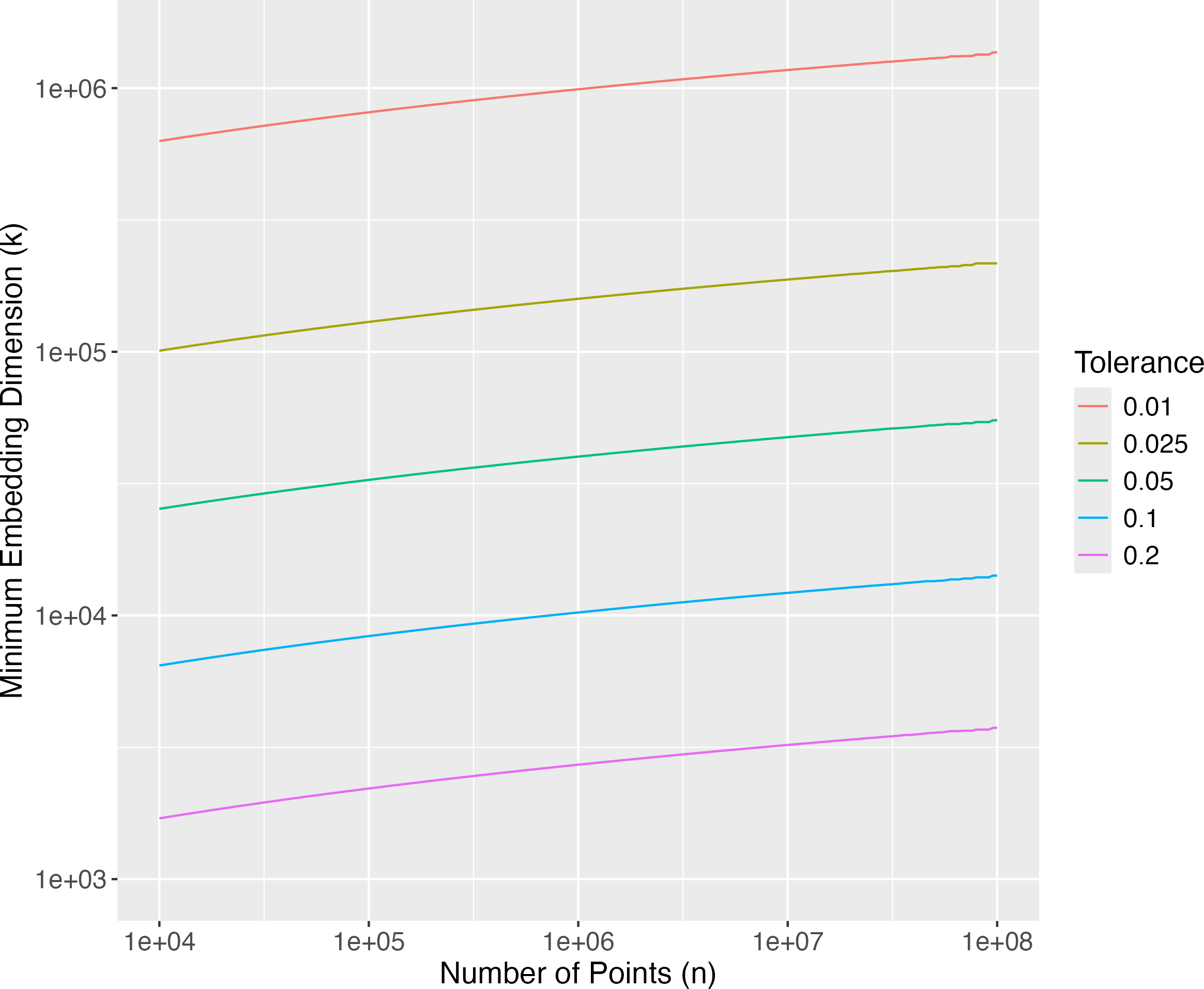}
\centering
\caption{The smallest embedding dimensions that lead to a positive no-failure probability using the bivariate or marginal lower bound.}
\label{fig:embedding_dimension}
\end{figure}

A practical question is whether the smallest embedding dimension that leads to a positive no-failure probability is greater for the bivariate or marginal lower bounds.
The bivariate lower bound is positive when
\begin{equation}
\mu>\sqrt{1-2\binom{n}{2}^{-1}},
\end{equation}
and the marginal lower bound is positive when
\begin{equation}
\label{mlb2}
\mu>1-\binom{n}{2}^{-1},
\end{equation}
so the goal is to compare the smallest values of $k$ that satisfy these inequalities.
To investigate this question, we compute the smallest embedding dimension that leads to a positive no-failure probability for five $\eps\in[0.01,0.2]$ and two-hundred $n\in[10^4,10^8]$.
The desired embedding dimensions are computed using bisection, and found to be the same for both bounds and all combinations of $n$ and $\eps$.
That is, the bivariate lower bound does not lead to a smaller embedding dimension than the marginal lower bound, which is reasonable considering the small difference between the bounds.
Figure~\ref{fig:embedding_dimension} shows the computed embedding dimensions.

In contrast to this numerical approach, the embedding dimension is often selected using an analytical lower bound derived with Markov's inequality \cite[eqn. 2.1]{dasgupta2003},
\begin{equation}
\label{dasgupta}
k\geq\frac{24\log(n)}{3\eps^2-2\eps^3}.
\end{equation}
Previously, \cite{li2006b} and \cite{rojo2010} found that embedding dimensions found numerically can be approximately 13-40\% lower than those found with analytical bounds.
For the values of $n$ and $\eps$ considered in this section, the value of $k$ found numerically is approximately 80-92\% of the value given by \eqref{dasgupta}.
\cite{li2006b} and \cite{rojo2010} appear to have considered different values of $n$ and $\eps$, which may account for the difference in these percentages.

\subsection{Data-dependent Lower Bound on the No-Failure Probability}
\label{subsec:data_bound}

Random projections are used for dimension reduction in part because they are independent of the data, unlike matrix factorization methods like singular value decomposition, for example.
Furthermore, their analysis is also usually independent of the data in that the embedding dimension and lower bounds on the no-failure probability do not depend on the data.
Both the marginal lower bound and the bivariate lower bound derived in Section \ref{sec:nfp} are independent of the data being embedded, for example.
However, a possibly greater data-dependent lower bound than the bivariate lower bound is given by \eqref{DeMorgan2},
\begin{equation}
\label{data_dependent_lower_bound}
p(\cap_{i=1}^{\binom{n}{2}}S_i)\geq 1 - \sum_{i=1}^{\frac{1}{2}\binom{n}{2}}\left(1-p(S_{2i-1}\cap S_{2i})\right),
\end{equation}
where the embedding dimension can be selected to ensure that the lower bound is positive.
A similar data-dependent lower bound on the no-failure probability is not available from the derivation of the marginal lower bound.
Actually computing this lower bound requires knowing the $\rho_{ii'}$, and so may be impossible if the number of points to embed or the original dimension of the data is large.
Also note that the indices of the lower bound can be permuted without changing the no-failure probability, so that the maximum can be taken over permutations of the indices.

\subsection{Limitations}
\label{subsec:limitations}

We briefly discuss two limitations of our result that also help situate this work in the broader context of random projection.
First, the lower bound on the no-failure probability relied on the assumption that all of the elements of the projection matrix were Gaussian distributed.
However, sparse, non-Gaussian random projections are often used since they can be implemented efficiently and approximately preserve distances \citep{achlioptas2003,li2006a}.
The bivariate lower bound derived here on the no-failure probability does not extend to these cases where pairs of projection errors do not have the bivariate gamma distribution.
Furthermore, the improvement in the lower bound on the no-failure probability found here is so small that it does not support using a Gaussian projection instead of a sparse projection.

A second limitation is that we derive a lower bound on the no-failure probability using the exact value of the failure probability instead of a bound on the failure probability.
This difference complicates a comparison of our work with the work of others on the JL lemma since using bounds obtained with Markov's inequality appears more common \citep{dasgupta2003}.
However, as suggested in \citet[footnote 7]{li2006b} and \cite{rojo2010}, the failure probability can be evaluated numerically and so does not need to be bounded using Markov's inequality.
One advantage with having an algebraic expression for the bound on the failure probability is that inverting the bound yields a closed-form expression for the embedding dimension.
This advantage is not helpful here though since our main goal is to bound the no-failure probability for a fixed embedding dimension, and not to determine an embedding dimension.

\section{Conclusion}
\label{sec:conclusion}

This paper considered the application of the bivariate gamma distribution to the proof of the JL lemma.
It was found that modifying the standard proof of the JL lemma to incorporate this distribution led to a small improvement in the lower bound on the no-failure probability.
Specifically, we followed the standard approach of using Bonferroni's inequality to bound the probability of at least one failure, where a failure is a significantly distorted distance.
We adapted this approach by applying Bonferroni's inequality to pairs of failures instead of single failures.
An inequality from \cite{jensen1969} that applies to the bivariate gamma distribution of two projection errors was then used to bound the probability of two successes.
One surprising aspect of this result is that it holds even if all of the pairwise differences between the points to be embedded are orthogonal.
Further work is needed to see if the result can be extended to higher-dimensional gamma distributions, or if the bivariate lower bound obtained here is the best lower bound that can be found with our approach.

\newpage

\appendix
\section{Appendix}
\label{app:appendix}

This appendix provides additional details on the bivariate gamma distribution of a pair of projection errors.
Section \ref{subsec:charfun} gives a derivation of the joint characteristic function of $V_i$ and $V_{i'}$, which is used to derive their covariance.
Section \ref{subsec:clt} then obtains a Gaussian approximation of the bivariate gamma distribution with the central limit theorem.
This approximation helps evaluate the bivariate gamma distribution function when the embedding dimension is large.
Section \ref{subsec:joint_dist} gives the joint distribution of all of the random projections and the conditional distribution of two random projections given a third.

\subsection{Characteristic Function of the Bivariate Gamma Distribution}
\label{subsec:charfun}

We now provide a derivation of the characteristic function of the bivariate gamma distribution and then use the function to derive the covariance of the distribution.
The characteristic or moment generating function is also given in \citet[eqn. 20]{wicksell1933}, \cite{kibble1941}, \cite{moran1967}, and \citet[eqn. 48.12]{kotz2004}, for example.
The more general characteristic or moment generating function of the multivariate gamma distribution can be found in \citet[eqn. 2.3]{krishnamoorthy1951} and \citet[eqn. 2]{krishnaiah1961}.

The joint characteristic function of $V_i$ and $V_{i'}$ can be derived from \eqref{gaussian_V}.
Since $Z_{i\ell}$ and $Z_{i'\ell'}$ are independent if $\ell\neq \ell'$ and the pairs $(Z_{i\ell},Z_{i'\ell})$ are identically distributed, the joint characteristic function satisfies
\begin{equation}
\phi_{ii'}(t_1,t_2)=[\theta_{ii'}(t_1,t_2)]^k,
\end{equation}
where $\theta_{ii'}(t_1,t_2)$ is the joint characteristic function of $Z_{ii'}=(Z_{i1},Z_{i'1})'$ and $t_1,t_2\in\mathbb{R}$.
The joint characteristic function of $Z_{ii'}$ is derived using standard algebraic manipulation and the `integrate-to-one' property of the Gaussian distribution.
Letting $\iota=\sqrt{-1}$ and
\begin{equation}
D=\begin{pmatrix}
t_1 & 0 \\
0 & t_2
\end{pmatrix},
\end{equation}
it follows that
\begin{align}
\theta_{ii'}(t_1,t_2)
&=\expval[\exp(\iota(t_1Z_{i1}^2+t_2Z_{i'1}^2))]\\
&=\expval[\exp(\iota Z_{ii'}'DZ_{ii'})]\\
&=\frac{1}{2\pi|\Sigma_{ii'}|^{1/2}}\int \exp\left(-\frac{1}{2}z_{ii'}'\Sigma_{ii'}^{-1}z_{ii'}+\iota z_{ii'}'Dz_{ii'}\right)\,dz_{ii'}\\
&=\frac{1}{2\pi|\Sigma_{ii'}|^{1/2}}\int \exp\left(-\frac{1}{2}z_{ii'}'\left[\Sigma_{ii'}^{-1}-2\iota D\right]z_{ii'}\right)\,dz_{ii'}\\
&=\frac{|R_{ii'}|^{1/2}}{|\Sigma_{ii'}|^{1/2}}\left(\frac{1}{2\pi|R_{ii'}|^{1/2}}\int \exp\left(-\frac{1}{2}z_{ii'}'R_{ii'}^{-1}z_{ii'}\right)\,dz_{ii'}\right),\qquad R_{ii'}^{-1}=\Sigma_{ii'}^{-1}-2\iota D\\
\label{kp1951}
&=|I_2-2\iota D\Sigma_{ii'}|^{-1/2}\cdot 1\\
&=((1-2\iota t_1)(1-2\iota t_2)+4t_1t_2\rho_{ii'}^2)^{-1/2},
\end{align}
where the one in \eqref{kp1951} is due to a bivariate normal density function being integrated over its support.

The joint characteristic function $\phi_{ii'}(t_1,t_2)$ is given in \citet[eqn. 6]{jensen1969} as an expression involving a matrix determinant,
\begin{align}
\phi_{ii'}(t_1,t_2)
&=
\begin{vmatrix}
(1-2\iota t_1)I_k & -2\iota t_1\rho_{ii'} I_k \\
-2\iota t_2\rho_{ii'} I_k & (1-2\iota t_2)I_k
\end{vmatrix}
^{-1/2},
\end{align}
where a scaling by 2 has been applied since the joint characteristic function of $V_i/2$ and $V_{i'}/2$ was given originally.
An equivalent expression for the joint characteristic function is
\begin{align}
\phi_{ii'}(t_1,t_2)
&=\left|(I_2-2\iota D\Sigma_{ii'})\otimes I_k\right|^{-1/2} \\
&=\left|I_2-2\iota D\Sigma_{ii'}\right|^{-k/2}|I_k|^{-2/2} \\
\label{chf}
&=((1-2\iota t_1)(1-2\iota t_2) + 4t_1t_2\rho_{ii'}^2)^{-k/2},
\end{align}
which is equivalent to $[\theta_{ii'}(t_1,t_2)]^k$, as expected.
The second equality is an application of an identity for the determinant of a Kronecker product.

The covariance of $V_i$ and $V_{i'}$ is given in \eqref{cov_proj_errors} and derived as follows.
First, $\expval[V_iV_{i'}]$ is computed from the characteristic function of $V_i$ and $V_{i'}$,
\begin{align}
\expval[V_iV_{i'}]
&=\left.-\frac{d}{dt_1}\frac{d}{dt_2}\phi(t_1,t_2)\right|_{t_1=t_2=0}\\
&=k^2+2k\rho_{ii'}^2.
\end{align}
Then, since $\expval[V_i]=\expval[V_{i'}]=k$, it follows that the covariance of $V_i$ and $V_{i'}$ is
\begin{align}
\cov[V_i,V_{i'}]
&=\expval[V_iV_{i'}]-\expval[V_i]\expval[V_{i'}]\\
&=2k\rho_{ii'}^2.
\end{align}
The result is useful since it shows how the data, through the dot product $\rho_{ii'}$, appears in the joint distribution of two projection errors.
Taking $k=1$ also implies that $\cov[Z_{i\ell}^2,Z_{i'\ell}^2]=2\rho_{ii'}^2$, as noted in Section \ref{subsec:bivgamma}.

\subsection{Gaussian Approximation of the Bivariate Gamma Distribution}
\label{subsec:clt}

We now detail the Gaussian approximation to the bivariate gamma distribution that is justified using the central limit theorem.
The basic idea is to approximate the asymptotic distribution of the right side of \eqref{gaussian_V} using the multivariate central limit theorem.
Following Theorem 1.17 in \cite{dasgupta2008}, the multivariate central limit theorem asserts that
\begin{equation}
\sqrt{k}\left(\frac{1}{k}
\begin{pmatrix}
V_i \\ V_{i'}
\end{pmatrix}
-
\begin{pmatrix}
1 \\ 1
\end{pmatrix}
\right)
\stackrel{d}{\rightarrow}
N(0_2,\Omega_{ii'})
\end{equation}
as $k\rightarrow\infty$, where the convergence is in distribution and the asymptotic covariance matrix is the covariance matrix of $(Z_{i1}^2,Z_{i'1}^2)'$,
\begin{equation}
\Omega_{ii'}=
\begin{pmatrix}
2 & 2\rho_{ii'}^2 \\
2\rho_{ii'}^2 & 2
\end{pmatrix}.
\end{equation}
Hence, if $k$ is large, then the joint distribution of $V_i$ and $V_{i'}$ can be approximated as
\begin{equation}
\begin{pmatrix}
V_i \\ V_{i'}
\end{pmatrix}
\stackrel{\cdot}{\sim}
N_2\left(\begin{pmatrix}
k \\ k
\end{pmatrix},
k\Omega_{ii'}
\right).
\end{equation}
Joint success probabilities can thus be computed by replacing the actual bivariate gamma distribution, $H(v_i,v_{i'})$, with this approximate Gaussian distribution.
For the results shown in this paper, we used the Gaussian approximation if $k>300$ and evaluated the joint distribution function directly using numerical integration if $k\leq300$.
This cut-off value of $k$ was justified by comparing values of the joint success probability evaluated using the Gaussian approximation and integrating \eqref{pdf_kibble} directly.
When $k=300$, the difference in the two values is less than one in ten-thousand for values of the correlation parameter across $(-1,1)$ and $\eps\in\{0.01, 0.025, 0.05\}$.

\subsection{Joint Distribution of Random Projections}
\label{subsec:joint_dist}

The $\binom{n}{2}$ random projections are jointly Gaussian distributed,
\begin{equation}
\label{full_joint_distn}
\begin{pmatrix}
    Z_1\\
    \vdots\\
    Z_{\binom{n}{2}}
  \end{pmatrix}
  \sim
  N(0_{\binom{n}{2}k},\Sigma\otimes I_k),
\end{equation}
where $\Sigma$ is a correlation matrix of dimension $\binom{n}{2}\times\binom{n}{2}$ with $(i,i')$th element $\rho_{ii'}$.
Since $\rho_{ii}=1$, the correlation matrix $\Sigma_{ii'}$ defined in \eqref{Sigma_ii} is the $2\times2$ submatrix formed by the $i$ and $i'$th rows and columns of $\Sigma$.

It follows from \eqref{full_joint_distn} that the joint distribution of three random projections, say $Z_1$, $Z_2$, and $Z_3$, is multivariate Gaussian,
\begin{equation}
\begin{pmatrix}
    Z_1\\
    Z_2\\
    Z_3
  \end{pmatrix}
  \sim
  N\left(0_{3k},
\begin{pmatrix}
1 & \rho_{12} & \rho_{13} \\
\rho_{12} & 1 & \rho_{23} \\
\rho_{13} & \rho_{23} & 1 \\
\end{pmatrix}
\otimes I_k\right).
\end{equation}
This implies that the conditional distribution of $Z_1$ and $Z_2$ given $Z_3$ is also multivariate Gaussian,
\begin{equation}
  \begin{pmatrix}
    Z_1\\
    Z_2
  \end{pmatrix}
   |Z_3\\
  \sim
  N\left(
  \begin{pmatrix}
    \rho_{13}\\
    \rho_{23}
  \end{pmatrix}\otimes Z_3,
\Sigma_{12\cdot3}\otimes I_k\right),
\end{equation}
where $\Sigma_{12\cdot3}$ is the covariance matrix of $Z_{1\ell}$ and $Z_{2\ell}$ given $Z_{3\ell}$ and given by
\begin{equation}
\Sigma_{12\cdot3}=
\begin{pmatrix}
1 - \rho_{13}^2 & \rho_{12}-\rho_{13}\rho_{23} \\
\rho_{12}-\rho_{13}\rho_{23} & 1 - \rho_{23}^2\\
\end{pmatrix}.
\end{equation}
Note that $\Sigma_{12\cdot3}$ is only a correlation matrix if the variances equal one.
The conditional correlation is thus
\begin{equation}
\label{partial_corr_coef}
\rho_{12\cdot3}=\frac{\rho_{12}-\rho_{13}\rho_{23}}{(1 - \rho_{13}^2)^{1/2}(1 - \rho_{23}^2)^{1/2}},
\end{equation}
which is called a partial correlation coefficient in the statistics literature \cite[sec. 4.5]{rencher2008}.
Hence, $Z_1$ and $Z_2$ are conditionally independent given $Z_3$ if and only if $\rho_{12\cdot3}=0$.

\section*{Acknowledgments}

\begin{quote}
The authors thank Geoffrey Sanders for helpful conversations on this topic.
\end{quote}

\section*{Disclaimer}

\begin{quote}
This document was prepared as an account of work sponsored by an agency  of the United States government. Neither the United States government nor Lawrence Livermore National Security, LLC, nor any of their employees makes any warranty, expressed or implied, or assumes any legal liability or responsibility for the accuracy, completeness, or usefulness of any information, apparatus, product, or process disclosed, or represents that its use would not infringe privately owned rights. Reference herein to any specific commercial product, process, or service by trade name, trademark, manufacturer, or otherwise does not necessarily constitute or imply its endorsement, recommendation, or favoring by the United States government or Lawrence Livermore National Security, LLC. The views and opinions of authors expressed herein do not necessarily state or reflect those of the United States government or Lawrence Livermore National Security, LLC, and shall not be used for advertising or product endorsement purposes.
\end{quote}

\begin{quote}
Lawrence Livermore National Laboratory is operated by Lawrence Livermore  National Security, LLC, for the U.S. Department of Energy, National Nuclear Security Administration under Contract DE-AC52-07NA27344.
\end{quote}

\clearpage

\bibliographystyle{apalike}
\bibliography{ms}

\end{document}